\documentclass[12pt,reqno]{article}

\usepackage[usenames]{color}
\usepackage[colorlinks=true,
linkcolor=webgreen, filecolor=webbrown,
citecolor=webgreen]{hyperref}

\definecolor{webgreen}{rgb}{0,.5,0}
\definecolor{webbrown}{rgb}{.6,0,0}

\usepackage{amssymb}
\usepackage{mathtools} 
\usepackage{graphicx}
\usepackage{amscd}
\usepackage{lscape}
\usepackage{tikz}
\usepackage{tikz-cd}
\usepackage{pgfplots}

\usetikzlibrary{matrix}
\usetikzlibrary{fit,shapes}
\usetikzlibrary{positioning, calc}
\tikzset{circle node/.style = {circle,inner sep=1pt,draw, fill=white},
        X node/.style = {fill=white, inner sep=1pt},
        dot node/.style = {circle, draw, inner sep=5pt}
        }
\usepackage{tkz-fct}

\usepackage{amsthm}
\newtheorem{theorem}{Theorem}

\theoremstyle{definition}

\newtheorem{example}[theorem]{Example}

\usepackage{float}

\usepackage{graphics,amsmath}
\usepackage{amsfonts}
\usepackage{latexsym}
\usepackage{epsf}

\newcommand{\seqnum}[1]{\href{http://oeis.org/#1}{\underline{#1}}}

\begin{document}

\begin{center}
\vskip 1cm{\LARGE\bf The Triple Riordan Group} \vskip 1cm \large
Paul Barry\\
School of Science\\
South East Technological University\\
Ireland\\
\href{mailto:pbarry@wit.ie}{\tt pbarry@wit.ie}
\end{center}
\vskip .2 in

\begin{abstract} We define the triple Riordan group, whose elements consist of  $4$-tuples of power series $(g, f_1, f_2, f_3)$ with $g\in \mathbf{R}[[x^3]]$, and
$f_1, f_2, f_3 \in x\mathbf{R}[[x^3]]$, for an appropriate ring $\mathbf{R}$. The construction of this group generalizes that of the double Riordan group, and lays the pattern for further generalizations.
\end{abstract}

The double Riordan group \cite{double, He}, which is a generalization of the checkerboard subgroup of the Riordan group \cite{book1, book2, SGWW}, has as its elements $3$-tuples of power series $(g,f_1,f_2)$ where
$g \in \mathbf{R}[[x^2]]$, $g_0 \ne 0$, and $f_1, f_2 \in x\mathbf{R}[[x^2]]$, $(f_i)_1 \ne 0$, where $\mathbf{R}$ is an appropriate ring, often taken to be $\mathbb{Z}$, or the fields $\mathbb{R}$ or $\mathbb{C}$. The element $(g, f_1, f_2)$ of the double Riordan group is then represented by the lower-triangular matrix $\left(a_{n,k}\right)$ whose $(n,k)$-th element is given by
$$ a_{n,k}=[x^n] g(x)f_1(x)^{\lfloor \frac{k+1}{2} \rfloor} f_2(x)^{\lfloor \frac{k}{2} \rfloor}.$$
Thus the columns of this matrix have generating functions $g, gf_1, gf_1f_2, gf_1^2f_2, gf_1^2f_2^2,\ldots$.

This approach can be generalized  to define the \emph{triple Riordan group} \cite{book2}. The (proper) triple Riordan group $\mathcal{T}\mathcal{R}$ is the group of $4$-tuples of power series $(g, f_1, f_2, f_3)$ with $g\in \mathbf{R}[[x^3]]$, $g_0 \ne 0$, and
$f_1, f_2, f_3 \in x\mathbf{R}[[x^3]]$, $(f_i)_1 \ne 0$, for an appropriate ring $\mathbf{R}$, with the following group operations. We let $h(x)=(f_1(x)f_2(x)f_3(x))^{\frac{1}{3}}$. Then the product in the group is given by
$$ \left(g, f_1, f_2, f_3\right)\cdot \left(G, F_1, F_2, F_3\right)=\left(gG(h), \frac{f_1}{h}F_1(h), \frac{f_2}{h}F_2(h), \frac{f_3}{h}F_3(h)\right).$$
For this product, we have following inverse of an element of the group.
$$\left(g, f_1, f_2, f_3\right)^{-1}=\left(\frac{1}{g(\bar{h})}, \frac{x\bar{h}}{f_1(\bar{h})}, \frac{x\bar{h}}{f_2(\bar{h})},\frac{x\bar{h}}{f_3(\bar{h})}\right).$$
Here, $\bar{h}$ represents the compositional inverse of $h$. The identity for the group is the element $(1,x,x,x)$.

The matrix that represents the element $\left(g,f_1, f_2, f_3\right)$ is the lower-triangular matrix $\left(a_{n,k}\right)$ with
$$a_{n,k}=[x^n] g(x)f_1(x)^{\lfloor \frac{k+2}{3} \rfloor}f_2(x)^{\lfloor \frac{k+1}{3} \rfloor}f_3(x)^{\lfloor \frac{k}{3} \rfloor}.$$

As seen in the product rule, the ``fundamental theorem of the triple Riordan group'' is that we have
$$\left(g,f_1,f_2,f_3\right)G=gG(h).$$

The bivariate generating function of the group element $\left(g,f_1, f_2,f_3\right)$ is given by
$$\frac{g(1+yf_1+y^2f_1f_2)}{1-y^3 f_1f_2f_3}.$$ In particular, the row sums and the diagonal sums of the associated matrix have generating functions
$$\frac{g(1+f_1+f_1f_2)}{1- f_1f_2f_3}$$ and
$$\frac{g(1+xf_1+x^2f_1f_2)}{1-x^3 f_1f_2f_3}$$ respectively.

We can identify some subgroups of the triple Riordan group $\mathcal{T}\mathcal{R}$. We have
\begin{align*}
\mathcal{A}&=\{(1,f_1,f_2,f_3) \,|\, f_i \in x\mathbf{R}[[x^3]], (f_i)_1 \ne 0\}\\
\mathcal{B}_1&=\{(g,xg,f_2,f_3)\,|\, g\in \mathbf{R}[[x^3]], f_i \in x\mathbf{R}[[x^3]], g_0 \ne 0, (f_i)_1 \ne 0\}\\
\mathcal{B}_2&=\{(g,f_1,xg,f_3)\,|\, g\in \mathbf{R}[[x^3]], f_i \in x\mathbf{R}[[x^3]], g_0 \ne 0, (f_i)_1 \ne 0\}\\
\mathcal{B}_3&=\{(g,xg,f_2,f_3)\,|\, g\in \mathbf{R}[[x^3]], f_i \in x\mathbf{R}[[x^3]], g_0 \ne 0, (f_i)_1 \ne 0\}.
\end{align*}
For instance, we have
\begin{align*}
(g,xg,f_2,f_3)\cdot (G, xG, F_2, F_3)&=\left(gG(h), \frac{xg}{h}(xG)(h), \frac{f_2}{h}F_2(h), \frac{f_3}{h}F_3(h)\right)\\
&=\left(gG(h), \frac{xg}{h}hG(h), \frac{f_2}{h}F_2(h), \frac{f_3}{h}F_3(h)\right)\\
&=\left(gG(h), xgG(h),  \frac{f_2}{h}F_2(h), \frac{f_3}{h}F_3(h)\right).\end{align*}
For the element $(g,x,x,x)$, we have $h(x)=x$, and so we have
$$(g,x,x,x)(1,f_1,f_2,f_3)=\left(g 1(x), \frac{x}{x}f_1(x), \frac{x}{x}f_2(x), \frac{x}{x}f_3(x)\right)=(g,f_1,f_2,f_3).$$
This means that the triple group $\mathcal{T}\mathcal{R}$ is the semi-direct product of the subgroup
$$\{ (g,x,x,x) \,|\, g \in \mathbf{R}[[x^3]], g_0 \ne 0\}$$ and $\mathcal{A}$.

Note that we have assumed above that $g_0=1, (f_i)_1=1, i=1,2,3$. This defines the \emph{proper} triple Riordan group, which, for $\mathbf{R}=\mathbb{Z}$, ensures integer entries for elements and their inverses. The following two examples follow this protocol. Relaxing these stipulations does not affect the group structure.
\begin{example}
We take the case of
\begin{align*}
g(x)&=\frac{1}{1-x^3}\\
f_1(x)&=\frac{x}{1-x^3}\\
f_2(x)&=x(1+x^3)\\
f_3(x)&=\frac{x}{1+x^3}.\end{align*}
This matrix begins
$$\left(\begin{array}{rrrrrrrrrr}
1 & 0 & 0 & 0 & 0 &  0 & 0 & 0 & 0 & \cdots \\
0 & 1 & 0 & 0 & 0 & 0 & 0 & 0 & 0 & \cdots \\
0 & 0 & 1 & 0 & 0 & 0 & 0 & 0 & 0 & \cdots \\
1 & 0 & 0 & 1 & 0 & 0  & 0 & 0 & 0 & \cdots \\
0 & 2 & 0 & 0 & 1 & 0  & 0 & 0 & 0 & \cdots \\
0 & 0 & 3 & 0 & 0 & 1 & 0 & 0 & 0 & \cdots \\
1 & 0 & 0 & 2 & 0 & 0  & 1 & 0 & 0 & \cdots \\
0 & 3 & 0 & 0 & 3 & 0  & 0 & 1 & 0 & \cdots \\
0 & 0 & 5 & 0 & 0 & 4 & 0 & 0 & 1 & \cdots \\
\vdots & \vdots & \vdots & \vdots & \vdots & \vdots & \vdots & \vdots & \vdots & \ddots
\end{array}\right).$$
The row sums of this matrix will have their generating function given by
$$s(x)=\frac{1+x+x^2-x^3+x^5}{(1-x^3)(1-2x^3)}.$$
This sequence, which begins
$$1, 1, 1, 2, 3, 4, 4, 7, 10, 8, 15, 22, 16, 31, 46, 32, 63, 94, 64, 127, 190,\ldots,$$
is seen to be an interleaving of $2^n$, $2^{n+1}-1$, and $3\cdot 2^n-2$ (\seqnum{A033484}).  This follows since
$$s(x)=\frac{1}{1-2x^3}+\frac{x}{(1-x^3)(1-2x^3)} +\frac{x^2(1+x^3)}{1-3x^3+2x^6}.$$ Note that we refer to known integer sequences by their entry number in the On-Line Encyclopedia of Integer Sequences (OEIS) \cite{SL1, SL2}.
We can illustrate the use of the fundamental theorem of triple Riordan arrays by looking at
$$\left(g, f_1, f_2, f_3\right)\cdot \frac{1-x^3}{1+x^3}.$$
We have
$$h(x)=(f_1(x)f_2(x)f_3(x))^{1/3}=\left(\frac{x^3}{1-x^3}\right)^{1/3}.$$
Then
$$\left(g, f_1, f_2, f_3\right)\cdot \frac{1-x^3}{1+x^3}=g(x)\frac{1-h(x)^3}{1+h(x)^3}=\frac{1-2x^3}{1-x^3},$$ which expands to give the sequence that begins
$$1,0,0,-1,0,0,-1,0,0,-1,0,0,-1,0,0,\ldots.$$
It is interesting to note that the reversion of this sequence gives a double aeration of the sequence \seqnum{A243659}, which gives the number of Sylvester classes of $3$-packed words of degree $n$ \cite{Novelli}.

The inverse of the above matrix is given by
$$\left(\frac{1}{1+x^3}, \frac{x}{1+x^3}, \frac{x(1+x^3)}{1+2x^3}, \frac{x(1+2x^3)}{1+x^3}\right).$$
The corresponding matrix begins
$$\left(\begin{array}{rrrrrrrrrr}
1 & 0 & 0 & 0 & 0 &  0 & 0 & 0 & 0 & \cdots \\
0 & 1 & 0 & 0 & 0 & 0 & 0 & 0 & 0 & \cdots \\
0 & 0 & 1 & 0 & 0 & 0 & 0 & 0 & 0 & \cdots \\
-1 & 0 & 0 & 1 & 0 & 0  & 0 & 0 & 0 & \cdots \\
0 &- 2 & 0 & 0 & 1 & 0  & 0 & 0 & 0 & \cdots \\
0 & 0 & -3 & 0 & 0 & 1 & 0 & 0 & 0 & \cdots \\
1 & 0 & 0 & -2 & 0 & 0  & 1 & 0 & 0 & \cdots \\
0 & 3 & 0 & 0 & -3 & 0  & 0 & 1 & 0 & \cdots \\
0 & 0 & 7 & 0 & 0 & -4 & 0 & 0 & 1 & \cdots \\
\vdots & \vdots & \vdots & \vdots & \vdots & \vdots & \vdots & \vdots & \vdots & \ddots
\end{array}\right).$$
The row sums of this inverse matrix have their generating function given by
$$\frac{1+x+x^2+3x^3+2x^4+x^5+2x^6}{(1+x^3)(1+2x^3)}.$$ The row sums sequence begins
$$1, 1, 1, 0, -1, -2, 0, 1, 4, 0, -1, -8, 0, 1, 16, 0, -1, -32, 0, 1, 64, 0, -1, \ldots.$$
This sequence is the interleaving of the three sequences $0^n$, $(-1)^n$ and $(-2)^n$.
We can operate on the expansion of $\frac{1-x^3}{1+x^3}$ as follows. We have in this case that
$$h(x)=\left(\frac{x^3}{1+x^3}\right)^{1/3}.$$
We find that
$$\left(\frac{1}{1+x^3}, \frac{x}{1+x^3}, \frac{x(1+x^3)}{1+2x^3}, \frac{x(1+2x^3)}{1+x^3}\right) \cdot \frac{1-x^3}{1+x^3}=\frac{1}{(1+x^3)(1+2x^3)}.$$

The sequence with generating function $$\frac{1}{(1+x^3)(1+2x^3)},$$  expands as
$$1, 0, 0, -3, 0, 0, 7, 0, 0, -15, 0, 0, 31, 0, 0, -63, 0, 0, 127, 0, 0,\ldots.$$
\end{example}
\begin{example}
For a second example we take the case of
\begin{align*}
g(x)&=1+x^3\\
f_1(x)&=x(1+x^3)\\
f_2(x)&=\frac{x}{1-x^3}\\
f_3(x)&=x(1-x^3).\end{align*}
Then the matrix defined by $\left(g,f_1,f_2,f_3\right)$ begins
$$\left(\begin{array}{rrrrrrrrrr}
1 & 0 & 0 & 0 & 0 &  0 & 0 & 0 & 0 & \cdots \\
0 & 1 & 0 & 0 & 0 & 0 & 0 & 0 & 0 & \cdots \\
0 & 0 & 1 & 0 & 0 & 0 & 0 & 0 & 0 & \cdots \\
1 & 0 & 0 & 1 & 0 & 0  & 0 & 0 & 0 & \cdots \\
0 & 2 & 0 & 0 & 1 & 0  & 0 & 0 & 0 & \cdots \\
0 & 0 & 3 & 0 & 0 & 1 & 0 & 0 & 0 & \cdots \\
0 & 0 & 0 & 2 & 0 & 0  & 1 & 0 & 0 & \cdots \\
0 & 1 & 0 & 0 & 3 & 0  & 0 & 1 & 0 & \cdots \\
0 & 0 & 4 & 0 & 0 & 4 & 0 & 0 & 1 & \cdots \\
\vdots & \vdots & \vdots & \vdots & \vdots & \vdots & \vdots & \vdots & \vdots & \ddots
\end{array}\right).$$
The row sums of this matrix will have their generating function given by
$$\frac{(1+x^3)(1+x+x^2-x^3+x^5-x^7)}{1-2x^3+x^9}.$$ This expands to give the sequence that begins
$$1, 1, 1, 2, 3, 4, 3, 5, 9, 5, 8, 17, 8, 13, 30, 13, 21, 51, 21, 34, 85, 34, 55, \ldots.$$
This is an interleaving of the sequences $F_{n+1}$ (\seqnum{A001519}$(n+1)$), $F_{2n+2}$ (\seqnum{A001906}$(n+1)$) and $F_{n+5}-4$ (\seqnum{A157728}).
Operating on the sequence with generating function $\frac{1+x^3}{1-x^3}$ we obtain the sequence with generating function
$$\frac{(1+x^3)(1+x^3+x^6)}{1-x^3-x^6}.$$
This is the (double) aeration of the sequence $2 F_{n+2}-\binom{1}{n}$, which has generating function
$$\frac{(1+x)(1+x+x^2)}{1-x-x^2}.$$
In order to describe the inverse of this matrix, we let
$$c(x)=\frac{1-\sqrt{1-4x}}{2x}$$ be the generating function of the Catalan numbers $C_n=\frac{1}{n+1}\binom{2n}{n}$ \seqnum{A000108}. The inverse matrix is then given by
$$\left(c(-x^3), xc(-x^3), x(1-x^3c(-x^3)), \frac{x}{1-x^3c(-x^3)}\right).$$ This matrix begins
$$\left(\begin{array}{rrrrrrrrrrr}
1 & 0 & 0 & 0 & 0 &  0 & 0 & 0 & 0 & 0 &\cdots \\
0 & 1 & 0 & 0 & 0 & 0 & 0 & 0 & 0 & 0 & \cdots \\
0 & 0 & 1 & 0 & 0 & 0 & 0 & 0 & 0 & 0 &\cdots \\
-1 & 0 & 0 & 1 & 0 & 0  & 0 & 0 & 0 & 0 &  \cdots \\
0 & -2 & 0 & 0 & 1 & 0  & 0 & 0 & 0 & 0 & \cdots \\
0 & 0 & -3 & 0 & 0 & 1 & 0 & 0 & 0 & 0 &  \cdots \\
2 & 0 & 0 & -2 & 0 & 0  & 1 & 0 & 0 & 0 & \cdots \\
0 & 5 & 0 & 0 & -3 & 0  & 0 & 1 & 0 & 0 & \cdots \\
0 & 0 & 8 & 0 & 0 & -4 & 0 & 0 & 1 & 0 & \cdots \\
-5 & 0 & 0 & 5 & 0 & 0 & -3 & 0 & 0 &1 & \cdots \\
\vdots & \vdots & \vdots & \vdots & \vdots & \vdots & \vdots & \vdots & \vdots & \ddots
\end{array}\right).$$
The row sums of this matrix then have a generating function given by
$$\frac{1+2x-x^2+x^3+3x^4+2x^5-2x^7-(1+2x-x^2-x^3-x^4)\sqrt{1+4x^3}}{2x^5(2-x^3)}.$$
This can also be written as
$$\frac{1+3x+2x^2+2x^3+2x^4-x^6}{1+2x-x^2+x^3+3x^4+2x^5-2x^7}c\left(\frac{x^5(2+6x+4x^2+3x^3+x^4-2x^5-4x^6-2x^7+x^9)}{(1+2x-x^2+x^3+3x^4+2x^5-2x^7)^2}\right).$$
The row sums sequence begins
$$1, 1, 1, 0, -1, -2, 1, 3, 5, -2, -8, -14, 6, 24, 42, -18, -75, -132, 57, 243, 429,\ldots.$$
This is the interleaving of three sequences. The first is $(-1)^n$ times the non-zero Fine numbers (\seqnum{A000957}) with generating function
$$\frac{c(-x)}{1-xc(-x)}=\frac{1-\sqrt{1+4x}}{x(\sqrt{1+4x}-3)}.$$
The second sequence is $(-1)^n$ times \seqnum{A000958} (which counts the number of ordered rooted trees with $n$ edges having root of odd degree), with generating function
$$\frac{c(-x)}{1-x^2c(-x)^2}=\frac{\sqrt{1+4x}-1}{x(1-2x+\sqrt{1+4x})}.$$ The third sequence is $(-1)^n C_{n+1}$, with generating function $\frac{1-c(-x)}{x}$.
The Hankel transforms \cite{Hankel} of these sequences are, respectively, the all $1$'s sequence, $F_{2n+1}$, and again the all $1$'s sequence.

\end{example}
\begin{example} For $g(x)\in \mathbf{R}[[x^3]]$ and $f(x) \in x\mathbf{R}[[x^3]]$, we have $(g,f)=(g,f,f,f)$. That is, the element $(g,f,f,f)$ of $\mathcal{T}\mathcal{R}$ is in fact a standard Riordan array. The element $\left(\frac{1}{1+x^3}, \frac{x}{1+x^3}\right)$ is an example of such an element, as is its inverse. The inverse begins
$$\left(\begin{array}{rrrrrrrrrrr}
1 & 0 & 0 & 0 & 0 &  0 & 0 & 0 & 0 & 0 &\cdots \\
0 & 1 & 0 & 0 & 0 & 0 & 0 & 0 & 0 & 0 & \cdots \\
0 & 0 & 1 & 0 & 0 & 0 & 0 & 0 & 0 & 0 &\cdots \\
1 & 0 & 0 & 1 & 0 & 0  & 0 & 0 & 0 & 0 &  \cdots \\
0 & 2 & 0 & 0 & 1 & 0  & 0 & 0 & 0 & 0 & \cdots \\
0 & 0 & 3 & 0 & 0 & 1 & 0 & 0 & 0 & 0 &  \cdots \\
3 & 0 & 0 & 4 & 0 & 0  & 1 & 0 & 0 & 0 & \cdots \\
0 & 5 & 0 & 0 & 5 & 0  & 0 & 1 & 0 & 0 & \cdots \\
0 & 0 & 12 & 0 & 0 & 6 & 0 & 0 & 1 & 0 & \cdots \\
12 & 0 & 0 & 18 & 0 & 0 & 7 & 0 & 0 &1 & \cdots \\
\vdots & \vdots & \vdots & \vdots & \vdots & \vdots & \vdots & \vdots & \vdots & \ddots
\end{array}\right).$$ This is \seqnum{A111373}, which enumerates left-factors of lattice paths in the positive quadrant, starting at $(0,0)$, that do not rise above the line $y=x$, with step set
$\{(1,1), (1,-2)\}$. The elements of this matrix $t_{n,k}$ obey the recurrence
$$t_{n,k}=t_{n-1,k-1}+t_{n-1,k+2}$$ with $t_{n,k}=0$ if $n<0$, or $k<0$, or $k>n$, along with $t_{0,k}=0^k$. The row sums \seqnum{A126042}
$$1, 1, 1, 2, 3, 4, 8, 13, 19, 38, 64, 98, 196, 337, 531, 1062, 1851,\ldots$$
count these left-factors up to the line $x=n$, as well as the number of ordered trees with $n+1$ edges, having non-root nodes of outdegree $0$ or $3$ (observation by E. Munarini).
\end{example}

\section{A lattice path example with three-fold step sets}
We have seen that the column generating functions of a triple Riordan array follow the pattern $g, gf_1, gf_1f_2, gf_1f_2f_3, gf_1^2f_2f_3,\ldots$. In this example, we follow this pattern, but we allow for the power series $f_i \in x\mathbb{Z}[[x]]$, $i=1,2,3$, with $g \in \mathbb{Z}[[x]]$. The resulting array in this case will be neither a Riordan array nor a triple Riordan array, but it will still have an interesting combinatorial meaning.

For this, we consider lattice paths in the positive quadrant, that begin at $(0,0)$, that do not rise above the diagonal $y=x$, and whose steps are governed as follows, where $t_{n,k}$ enumerates the number of paths from $(0,0)$ to $(n,k)$. We have
\[
t_{n,k} =
\begin{dcases*}
 t_{n-1,k-1}+t_{n-1,k+1}
   & if  $k \bmod 3=0$ \\[1ex]
 t_{n-1,k-1}+t_{n-1,k}+t_{n-1,k+1}
   & if $k \bmod 3=1$\\[1ex]
 t_{n-1,k-1}+t_{n-2,k}+t_{n-1,k+1} 
 & if $k \bmod 3=2$.
\end{dcases*}
\]
Thus the step sets alternate between 
\begin{align*}&\text{Dyck steps}\quad &\{(1,1), (1,-1)\},\\ 
&\text{Motzkin steps}\quad &\{(1,1),(1,0),(1,-1)\}, \\
&\text{Schroeder steps}\quad &\{(1,1),(2,0),(1,-1)\}.\end{align*}  Boundary conditions are that 
$t_{n,k}=0$ for $n<0$ or $k<0$ or $k>n$, and $t_{n,0}=0^k$. The resulting matrix $\left(t_{n,k}\right)$ begins 
$$\left(\begin{array}{rrrrrrrrrrr}
1 & 0 & 0 & 0 & 0 &  0 & 0 & 0 & 0 & 0 &\cdots \\
0 & 1 & 0 & 0 & 0 & 0 & 0 & 0 & 0 & 0 & \cdots \\
1 & 1 & 1 & 0 & 0 & 0 & 0 & 0 & 0 & 0 &\cdots \\
1 & 3 & 1 & 1 & 0 & 0  & 0 & 0 & 0 & 0 &  \cdots \\
3 & 5 & 5 & 1 & 1 & 0  & 0 & 0 & 0 & 0 & \cdots \\
5 & 13 & 7 & 6 & 2 & 1 & 0 & 0 & 0 & 0 &  \cdots \\
13 & 25 & 24 & 9 & 9 & 2  & 1 & 0 & 0 & 0 & \cdots \\
25 & 62 & 41 & 33 & 20 & 11  & 2 & 1 & 0 & 0 & \cdots \\
62 & 128 & 119 & 61 & 64 & 24 & 12 & 3 & 1 & 0 & \cdots \\
128 & 309 & 230 & 183 & 149 & 87 & 27 & 16 & 3 &1 & \cdots \\
\vdots & \vdots & \vdots & \vdots & \vdots & \vdots & \vdots & \vdots & \vdots & \ddots
\end{array}\right).$$
Then we have 
$$t_{n,k}=[x^n]g(x)f_1(x)^{\lfloor \frac{k+2}{3} \rfloor}f_2(x)^{\lfloor \frac{k+1}{3} \rfloor}f_3(x)^{\lfloor \frac{k}{3} \rfloor},$$ 
where 
\begin{align*}
g(x)&=\frac{1-x-2x^2+x^3}{1-x-2x^2+x^4}c\left(\frac{x^2(1-x-x^2)(1-x-2x^2+x^3)}{(1-x-2x^2+x^4)^2}\right),\\
f_1(x)&=\frac{1}{x}\left(1-\frac{1}{g(x)}\right),\\
f_2(x)&=\frac{x(1-x-x^2)}{1-x-2x^2+2x^3+x^4}c\left(\frac{x^2(1-x-x^2)(1-2x^2)}{(1-x-2x^2+2x^3+x^4)^2}\right),\\
f_3(x)&=xg(x).\end{align*}
We conclude that the left factors of these paths to the line $x=n$ have generating function 
$$\frac{1-2x^2}{1-2x-3x^2+4x^3+x^4-x^5}c\left(\frac{-x(1-x-x^2)(1-2x^2)(1-x-4x^2+x^4)}{(1-2x-3x^2+4x^3+x^4-x^5)^2}\right).$$ 
This generating function expands to give the sequence that begins 
$$1, 1, 3, 6, 15, 34, 83, 195, 474, 1133, 2756,\ldots.$$

\section{The quadruple Riordan group}
Taking into account the definitions of the double and the triple Riordan groups, we can define the quadruple Riordan group to be the set of elements $(g,f_1,f_2,f_3,f_4)$ with
$g \in \mathbf{R}[[x^4]]$, and $f_i \in x\mathbf{R}[[x^4]]$, $i=1,2,3,4$,  where the lower-triangular matrix $(a_{n,k})$ that represents the element $(g,f_1,f_2,f_3,f_4)$ is given by
$$a_{n,k}=[x^n]g(x)f_1^{\lfloor \frac{k+3}{4} \rfloor}f_2^{\lfloor \frac{k+2}{4} \rfloor} f_3^{\lfloor \frac{k+1}{4} \rfloor}f_4^{\lfloor \frac{k}{4} \rfloor}.$$
We define $h(x)=\left(f_1(x)f_2(x)f_3(x)f_4(x)\right)^{1/4}$. Then the product rule for the group is given by
$$ \left(g, f_1, f_2, f_3, f_4\right)\cdot \left(G, F_1, F_2, F_3,F_4\right)=\left(gG(h), \frac{f_1}{h}F_1(h), \frac{f_2}{h}F_2(h), \frac{f_3}{h}F_3(h), \frac{f_4}{h}F_4(h)\right).$$
The inverse in the group is then given by
$$\left(g, f_1, f_2, f_3, f_4\right)^{-1}=\left(\frac{1}{g(\bar{h})}, \frac{x\bar{h}}{f_1(\bar{h})}, \frac{x\bar{h}}{f_2(\bar{h})},\frac{x\bar{h}}{f_3(\bar{h})},\frac{x\bar{h}}{f_4(\bar{h})}\right).$$
The identity in the group is given by $(1,x,x,x,x)$.

Higher order groups can be defined in a similar manner \cite{book2}.

\bigskip
\hrule
\bigskip
\noindent 2020 {\it Mathematics Subject Classification}:
Primary 15B36; Secondary 05A15, 11B83, 11C20, 15A15.

\noindent \emph{Keywords:} Riordan array, double Riordan group, triple Riordan group, generating function, integer sequence.

\bigskip
\hrule
\bigskip
\noindent (Concerned with sequences
\seqnum{A000108},
\seqnum{A000957},
\seqnum{A000958},
\seqnum{A001519},
\seqnum{A001906},
\seqnum{A033484},
\seqnum{A111373},
\seqnum{A126042},
\seqnum{A157728}, and
\seqnum{A243659}).

\end{document}